\documentclass[12pt]{article}
\usepackage{amsfonts,amssymb,latexsym}
\newtheorem{theorem}{Theorem}[section]
\newtheorem{lemma}[theorem]{Lemma}
\newtheorem{corollary}[theorem]{Corollary}
\newtheorem{proposition}[theorem]{Proposition}

\newenvironment{proof}
{\par\addvspace{0.3cm}\noindent{\rm Proof. }}
{\nopagebreak\mbox{}\hfill $\Box$\par\addvspace{0.25cm}}

\renewcommand{\Re}{\mbox{\rm Re\,}}
\renewcommand{\Im}{\mbox{\rm Im\,}}
\newcommand{\diag}{\mbox{\rm diag\,}}
\newcommand{\im}{\mbox{\rm im\,}}

\newcommand{\cR}{{\cal R}}
\newcommand{\R}{{\mathbb R}}
\newcommand{\Rp}{{\mathbb R}_+}
\newcommand{\C}{{\mathbb C}}
\newcommand{\Z}{{\mathbb Z}}
\newcommand{\T}{{\mathbb T}}
\renewcommand{\kappa}{\varkappa}
\newcommand{\qed}{\hfill $\Box$}
\newcommand{\be}{\begin{equation}}
\newcommand{\ee}{\end{equation}}

\newcommand{\bq}{\begin{eqnarray}}
\newcommand{\eq}{\end{eqnarray}}
\newcommand{\nn}{\nonumber}
\newcommand{\ba}{\begin{array}}
\newcommand{\ea}{\end{array}}

\newcommand{\wt}[1]{\widetilde{#1}}
\newcommand{\iv}{^{-1}}
\newcommand{\iy}{\infty}

\newcommand{\ovl}{\overline}
\newcommand{\Hp}{H^p(\T)}
\newcommand{\Hq}{H^q(\T)}
\newcommand{\Lp}{L^p(\T)}
\newcommand{\Lq}{L^q(\T)}
\newcommand{\Lpe}{L^p_{\rm even}(\T)}
\newcommand{\Lqe}{L^q_{\rm even}(\T)}
\newcommand{\Li}{L^\iy(\T)}
\newcommand{\cL}{{\cal L}}

\newcommand{\cM}{\mathcal{M}}
\newcommand{\cS}{\mathcal{S}}

\begin{document}

\date{}
\title{Factorization theory for a class of Toeplitz + Hankel operators}
\author{Estelle L. Basor\thanks{ebasor@calpoly.edu. 
          Supported in part by NSF Grant DMS-9970879.}\\
               Department of Mathematics\\
               California Polytechnic State University\\
               San Luis Obispo, CA 93407, USA
        \and
        Torsten Ehrhardt\thanks{tehrhard@mathematik.tu-chemnitz.de.}\\
       		Fakult\"{a}t f\"{u}r Mathematik\\
         	Technische Universit\"{a}t Chemnitz\\
         	09107 Chemnitz, Germany}
\maketitle

\begin{abstract}
In this paper we study operators of the form $M(\phi)=T(\phi)+H(\phi)$
where $T(\phi)$ and $H(\phi)$ are the Toeplitz and Hankel operators
acting on $H^p(\T)$ with generating function $\phi\in L^\iy(\T)$.
It turns out that $M(\phi)$ is invertible if and only if the function
$\phi$ admits a certain kind of generalized factorization. 
\end{abstract}

%%%%%%%%%%%%%%%%%%%%%%%%%%%%%%%%%%%%%%%%%%%%%%%%%%%%%%%%%%%%%%%%%%
%%%%%%%%%%%%%%%%%%%%%%%%%%%%%%%%%%%%%%%%%%%%%%%%%%%%%%%%%%%%%%%%%%
%%%%%%%%%%%%%%%%%%%%%%%%%%%%%%%%%%%%%%%%%%%%%%%%%%%%%%%%%%%%%%%%%%

\section{Introduction}
\noindent
This paper is devoted to the study of operators of the form
\bq M(\phi) &=& T(\phi)+H(\phi) \eq
acting on the Hardy space $H^p(\T)$ where $1<p<\iy$.
Here $\phi\in L^\iy(\T)$ is a Lebesgue measurable and essentially 
bounded function on the unit circle $\T$.
The Toeplitz and Hankel operators are defined by
\be\label{f.2}
T(\phi):f\mapsto P(\phi f),\qquad
H(\phi):f\mapsto P(\phi (Jf)),
\qquad f\in H^p(\T),
\ee
where is $J$ the following flip operator,
\be\label{f.3}
J:f(t)\mapsto t\iv f(t\iv),\quad t\in\T,
\ee
acting on the Lebesgue space $L^p(\T)$. 
The operator $P$ stands for the Riesz projection,
\be
P:\sum_{n=-\iy}^\iy f_nt^n\mapsto \sum_{n=0}^\iy f_nt^n, \quad t\in\T,
\ee
which is bounded on $L^p(\T)$, $1<p<\iy$, and whose image is $H^p(\T)$.
The complex conjugate Hardy space $\ovl{\Hp}$ is the set of all
functions $f$ whose complex conjugate belongs to $\Hp$.

For a Banach algebra $B$ (such as $L^\iy(\T)$ or $H^\iy(\T)$) we will denote
by $GB$ the group of all invertible elements.

It is well known that a necessary and sufficient condition for the
Fredholmness and the invertibility of the Toeplitz operator $T(\phi)$
with $\phi\in L^\iy(\T)$ can be given in terms of the Wiener-Hopf
factorability of the generating function $\phi$. 
Let us recall the underlying definitions and results
\cite{CG,GK,LS} (see also \cite{Wi} for the case $p=2$).

A function $\phi\in L^\iy(\T)$ is said to admit a
{\em factorization in $L^p(\T)$} if one can write
\bq
\phi(t) &=& \phi_-(t)t^\kappa\phi_+(t),\qquad t\in\T,
\eq
where $\kappa$ is an integer and the factors $\phi_-$ and $\phi_+$ satisfy the
following conditions:
\begin{itemize}
\item[(i)]
$\phi_-\in \ovl{\Hp}$, $\phi_-\iv\in \ovl{\Hq}$,
\item[(ii)]
$\phi_+\in\Hq$, $\phi_+\iv\in\Hp$,
\item[(iii)]
The linear operator $f\mapsto \phi_+\iv P(\phi_-\iv f)$, which
is defined on the set of all trigonometric polynomials and takes values in
$\Lp$, can be extended by continuity to a linear bounded operator
acting from $L^p(\T)$ into $L^p(\T)$.
\end{itemize}
Here $p\iv+q\iv=1$. 
A factorization where merely (i) and (ii) but not necessarily
condition (iii) is fulfilled is called a
{\em weak factorization in $L^p(\T)$}.
The number $\kappa$ is called the index of the (weak)
factorization and is uniquely determined.

The crucial result is that for given $\phi\in\Li$ the operator $T(\phi)$ is
a Fredholm operator on the space $\Hp$ if and only if the function
$\phi$ admits a factorization in $\Lp$. In this case the defect numbers are
given by
\be
\dim\ker T(\phi)=\max\{0,-\kappa\},\qquad
\dim\ker (T(\phi))^*=\max\{0,\kappa\}.
\ee
Hence $T(\phi)$ is invertible if and only if $\phi$ admits a factorization
in $\Lp$ with index $\kappa=0$.

The goal of the present paper is to obtain a corresponding results
for operators $M(\phi)$ with $\phi\in\Li$. We will encounter another type 
of factorization of the function $\phi$ which is related to the Fredholm
theory of these operators. The investigations taken up in
this paper are motivated by the results obtained in \cite{BaEh1,E}.

%%%%%%%%%%%%%%%%%%%%%%%%%%%%%%%%%%%%%%%%%%%%%%%%%%%%%%%%%%%%%%%%%%
%%%%%%%%%%%%%%%%%%%%%%%%%%%%%%%%%%%%%%%%%%%%%%%%%%%%%%%%%%%%%%%%%%
%%%%%%%%%%%%%%%%%%%%%%%%%%%%%%%%%%%%%%%%%%%%%%%%%%%%%%%%%%%%%%%%%%

\section{Basic properties of $M(\phi)$}

It is well known that for Toeplitz and Hankel operators the following relations 
hold:
\bq
T(\phi\psi) &=& 
T(\phi)T(\psi)+H(\phi)H(\wt{\psi}),\label{f2.Tx}\\[.2ex]
H(\phi\psi) &=& T(\phi)H(\psi)+H(\phi)T(\wt{\psi}).\label{f2.Hx}
\eq
Here $\wt{\psi}(t)=\psi(t\iv)$.
Adding both equations, it follows that
$$ M(\phi\psi)\;=\; T(\phi)M(\psi)+H(\phi)M(\wt{\psi}),$$
and hence
\bq\label{f.Mx}
M(\phi\psi) &=& M(\phi)M(\psi)+H(\phi)M(\wt{\psi}-\psi).\label{f2.Mx}
\eq

In some particular cases, this identity simplifies to a multiplicative relation:
\bq\label{f.M}
M(\phi\psi) &=& M(\phi)M(\psi)
\eq
if $\phi\in\ovl{H^\iy(\T)}$ or $\psi=\wt{\psi}$.
Based on this identity one can establish a sufficient invertibility criteria
for $M(\phi)$, which anticipates to some extent the factorization result
that we establish in this paper. 
Assume that $\phi\in\Li$ admits a factorization $\phi=\phi_-\phi_0$, where
$\phi_-\in G\ovl{H^\iy(\T)}$ and $\phi_0\in G\Li$ such that $\wt{\phi}_0=\phi_0$. 
Then $M(\phi)$ is invertible and its inverse its given by
$M(\phi_0\iv)M(\phi_-\iv)$.

It is occasionally convenient to introduce the multiplication operator
on $\Lp$,
\be
L(\phi):\Lp\to\Lp,f\mapsto \phi f, 
\ee
and then consider Toeplitz and Hankel operators as
restrictions of the following operators onto $\Hp$:
\be
T(\phi) = PL(\phi)P|_{\Hp},\qquad
H(\phi) = PL(\phi)JP|_{\Hp}.
\ee
The Toeplitz + Hankel operators $M(\phi)$ are of the form
\be\label{f.Mop}
M(\phi) = PL(\phi)(I+J)P|_{\Hp}.
\ee

The following results gives estimates for the norm of $M(\phi)$, where the 
constant $C_p$ depends only on the parameter $p$.

\begin{proposition}\label{p2.1}
Let $\phi\in L^\iy(\T)$. Then
$\;\|\phi\|_{\Li} \le\|M(\phi)\|_{\cL(\Hp)} \le C_p\|\phi\|_{\Li}$.
\end{proposition}
\begin{proof}
The upper estimate follows from (\ref{f.Mop}) Note that
the operators $P$ and $J$ are both bounded on $\Lp$.
In order to prove the lower estimate, put $U_n=M(t^n)$.
Since $U_nU_{-n}=I$ and $JU_n=U_{-n}J$, it is easy to see that
$$ U_{-n}M(\phi)U_n \;=\; (U_{-n}PU_n)L(\phi)(U_{-n}PU_n)+
(U_{-n}PU_{-n})L(\phi)J(U_{-n}PU_n).$$
Using the fact that $U_{-n}PU_n\to I$ and $U_{-n}PU_{-n}\to 0$ strongly
on $L^p(\T)$ as $n\to\iy$, it follows that
\bq U_{-n}M(\phi)U_n &\to& L(\phi) \label{f2.UMU}\eq
strongly on $L^p(\T)$ as $n\to\iy$. Since $U_{\pm n}$ are isometries on
$L^p(\T)$, this implies the lower estimate.
\end{proof}

Now we obtain a necessary condition for the Fredholmness of $M(\phi)$.

\begin{proposition}\label{p2.2}
Suppose that $M(\phi)$ is Fredholm on $\Hp$. Then
$\phi \in GL^\iy(\T)$.
\end{proposition}
\proof
The proof is based on standard arguments.
If $M(\phi)$ is a Fredholm operator, then there exist a
$\delta>0$  and a finite rank projection $K$ on the kernel of 
$M(\phi)$
such that
\bq
\|M(\phi)f\|_{H^p(\T)}+\|Kf\|_{H^p(\T)} &\ge& \delta\|f\|_{H^p(\T)}\nn
\eq
for all $f\in H^p(\T)$. Putting $Pf$ instead of $f$, this implies that
\bq
\|M(\phi)f\|_{L^p(\T)}+\|KPf\|_{L^p(\T)}+\delta\|(I-P)f\|_{L^p(\T)}
&\ge& \delta\|f\|_{L^p(\T)}\nn
\eq
for all $f\in L^p(\T)$. Replacing $f$ by $U_nf$ and observing again 
that
$U_{\pm n}$ are isometries on $L^p(\T)$, it follows that
\bq
\|U_{-n}M(\phi)U_nf\|_{L^p(\T)}+\|KPU_nf\|_{L^p(\T)}+\delta
\|U_{-n}(I-P)U_nf\|_{L^p(\T)} &\ge& \delta\|f\|_{L^p(\T)}.\nn
\eq
Now we take the limit $n\to\iy$. Because $U_n\to0$ weakly, we have 
$KPU_n\to0$
strongly. It remains to apply (\ref{f2.UMU}) and again the fact that
$U_{-n}PU_n\to I$ strongly. We obtain
\bq
\|L(\phi) f\|_{L^p(\T)} &\ge& \delta\|f\|_{L^p(\T)}.\nn
\eq
This implies that $\phi\in G\Li$.
\qed

The following lemma, which appears in slightly less general form as 
Lemma 2.9 in \cite{BS2}, turns out to be useful.

\begin{lemma}\label{l2.3}
Let $X_1$ and $X_2$ be linear spaces, $A:X_1\to X_2$ be a linear
and invertible operator, $P_1:X_1\to X_1$ and $P_2:X_2\to X_2$ be linear
projections, and $Q_1=I-P_1$ and $Q_2=I-P_2$. Then
$P_2AP_1:\im P_1\to \im P_2$ is invertible if and only if
$Q_1A\iv Q_2:\im Q_2\to \im Q_1$ is invertible.
\end{lemma}
\begin{proof}
The following formulas, which relate both inverses with each other,
can be verified by a direct calculation:
\bq
(P_2AP_1)\iv &=& P_1A\iv P_2- P_1A\iv Q_2(Q_1 A\iv Q_2)\iv Q_1A\iv P_2,\nn\\
(Q_1 A\iv Q_2)\iv &=& Q_2AQ_1-Q_2AP_1(P_2AP_1)\iv P_2AQ_1.\nn
\eq
The validity of these formulas implies the desired assertion.
\end{proof}

In what follows, we are going to relate $M(\phi)$ with two other
operators $\Phi(\phi)$ and $\Psi(\phi)$. First of all note that
$J^2=I$. Hence 
\be
P_{J}=\frac{I+J}{2},\qquad
Q_{J}=\frac{I-J}{2}
\ee
are complementary projections acting on the space $\Lp$ and decompose this
space into the direct sum $\Lp=\im P_J|_{\Lp}\dotplus\im Q_J|_{\Lp}$.
Let $L^p_J(\T)=\im P_J|_{\Lp}$. Given $\phi\in\Li$, define
\be\label{f.12}
\Phi(\phi)=PL(\phi)P_J,\qquad
\Psi(\phi)=P_JL(\phi)P,
\ee
where we consider these operators as acting between the following spaces:
\be
\Phi(\phi):L^p_J(\T)\to\Hp,\qquad
\Psi(\phi):\Hp\to L^p_J(\T),
\ee

\begin{proposition}\label{p2.4}
Let $\phi\in G\Li$. Then the following assertions are equivalent:
\begin{itemize}
\item[(i)]
$M(\phi)$ is invertible in $\cL(\Hp)$,
\item[(ii)]
$\Phi(\phi)$ is invertible in $\cL(L^p_J(\T),\Hp)$,
\item[(iii)]
$\Psi(\psi)$ is invertible in $\cL(\Hp,L^p_J(\T))$,
where $\psi(t)=\phi\iv(-t\iv)$.
\end{itemize}
\end{proposition}
\begin{proof}
In order to prove the equivalence of (i) and (ii) we refer to the
representation (\ref{f.Mop}) and (\ref{f.12}) of these operators.
Moreover, note that $A=(I+J)P:\Hp\to L^p_J(\T)$ and
$B=\frac{1}{2}P(I+J):L^p_J(\T)\to\Hp$ are inverse to each other:
\bq
AB &=& \frac{1}{2}(I+J)P(I+J)\;=\;\frac{P+JP+Q+JQ}{2}\;=\;\frac{I+J}{2}\nn\\
BA &=& \frac{1}{2}P(I+J)^2P\;=\; P(I+J)P\;=\;P\nn
\eq
Here we have used the basic relations, $J^2=I$ and $JPJ=Q$, where
$Q=I-P$.

In order to prove the equivalence of (ii) and (iii) we apply Lemma \ref{l2.3}
with $P_1=P_J$, $P_2=P$, $Q_1=Q_J$ and $Q_2=Q$. We obtain that
$\Phi(\phi)$ is invertible if and only if the operator
$$
Q_JL(\phi\iv)Q:\im Q|_{\Lp}\to\im Q_J|_{\Lp}
$$
is invertible. Multiplying this operator from the left and right with
$J$, we obtain the operator
$$
Q_JJL(\phi\iv)JP:\Hp\to\im Q_J|_{\Lp}
$$
since $QJ=JP$ and $JQ_J=Q_JJ$. Next we multiply the latter operator
from the left and right with $W$, where $W:f(t)\mapsto f(-t)$.
We obtain the operator
$$
P_JWJL(\phi\iv)JWP:\Hp\to L^p_J(\T)
$$
since $WQ_J=P_JW$ (i.e., $WJ+JW=0$) and $PW=WP$. This last operator is
equal to $\Psi(\psi)$ since $WJL(\phi\iv)JW=L(\psi)$ as one can readily
check. From this it follows that
$\Phi(\phi)$ is invertible if and only if so is $\Psi(\psi)$.
\end{proof}

%%%%%%%%%%%%%%%%%%%%%%%%%%%%%%%%%%%%%%%%%%%%%%%%%%%%%%%%%%%%%%%%%%
%%%%%%%%%%%%%%%%%%%%%%%%%%%%%%%%%%%%%%%%%%%%%%%%%%%%%%%%%%%%%%%%%%
%%%%%%%%%%%%%%%%%%%%%%%%%%%%%%%%%%%%%%%%%%%%%%%%%%%%%%%%%%%%%%%%%%

\section{Weak asymmetric and antisymmetric factorization}

In what follows we introduce the notion of a weak asymmetric
and a weak antisymmetric factorization of a matrix function
$\phi\in G\Li$ in the space $\Lp$.
Before we do so, however, let us for a 
moment try to anticipate the kind of factorization we should expect 
based on our notions from the theory of Toeplitz operators.
If we suppose that $M(\phi)$ is invertible on $\Hp$,
then there exists a function $h\in\Hp$ such that $M(\phi)h = 1$. Now this 
means, using the definition of $J$ that  
$$\phi(h+t^{-1}\tilde{h}) = g_{-}$$
where $g_{-} \in\ovl{\Hp}$.
Multiplying this equation with the factor $(1+t\iv)\iv$ gives
$$\phi(1+t\iv)\iv(h+t^{-1}\tilde{h}) = (1+t\iv)\iv g_{-}.$$
Introduce the functions
$$
f_0=(1+t\iv)\iv(h+t^{-1}\tilde{h})\quad\mbox{ and }\quad
f_-=(1+t\iv)\iv g_{-}.
$$
Thus $\phi f_0=f_-$, which gives upon assuming for a moment
the invertibility of $f_0$ the factorization $\phi=f_0\iv f_-$.
On the other hand
$$
\tilde{f}_0=(1+t)\iv(\tilde{h}+th)=(1+t\iv)\iv(h+t^{-1}\tilde{h})=f_0.
$$
Hence $f_0$ is even. From the definition of $f_-$ and $f_0$ it now follows that
$$
(1+t\iv)f_-\in\ovl{\Hp}\quad\mbox{ and }\quad
|1+t|\iv f_0\in \Lpe.
$$
Here and in what follows $\Lpe$ stands for the set of all functions 
$\phi_0\in\Lp$ which are even, i.e., $\phi_0=\tilde{\phi}_0$.

The above analysis shows that the factor $1+t^{-1}$ plays a special role in any 
factorization theory and hints of its presence in the following definitions.
It turns out that the same is true for the
factor $1-t^{-1}$. This is not apparent from the heuristics that 
have just been presented, but will be become clear later on.

A function $\phi\in G\Li$ is said to admit a {\em weak asymmetric
factorization in $\Lp$} if it admits a representation
\bq
\phi(t) &=&\phi_-(t)t^\kappa\phi_0(t),\qquad t\in\T,
\eq
such that $\kappa\in\Z$ and
\begin{itemize}
\item[(i)]
$(1+t\iv)\phi_-\in\ovl{\Hp}$, $(1-t\iv)\phi_-\iv\in\ovl{\Hq}$,
\item[(ii)]
$|1-t|\phi_0\in\Lqe$, $|1+t|\phi_0\iv\in\Lpe$,
\end{itemize}
where $p\iv+q\iv=1$.

The uniqueness of a weak asymmetric factorization (up to a constant)
is stated in the following proposition.

\begin{proposition}\label{p3.1}
Assume that $\phi$ admits two weak asymmetric factorizations in $\Lp$:
\be\label{f3.2}
\phi(t)\;=\;\phi_-^{(1)}(t)t^{\kappa_1}\phi_0^{(1)}(t)\;=\;
\phi_-^{(2)}(t)t^{\kappa_2}\phi_0^{(2)}(t),\qquad t\in\T.
\ee
Then $\kappa_1=\kappa_2$ and $\phi_-^{(1)}=\gamma\phi_-^{(2)}$,
$\phi_0^{(1)}=\gamma\iv\phi_0^{(2)}$ with $\gamma\in\C\setminus\{0\}$.
\end{proposition}
\begin{proof}
Assume without loss of generality that $\kappa=\kappa_1-\kappa_2\le0$. 
{}From (\ref{f3.2}) it follows that 
\bq\label{f3.16}
(\phi_-^{(2)})\iv\phi_-^{(1)}t^\kappa&=&\phi_0^{(2)}(\phi_0^{(1)})\iv.
\eq
Put $\psi=(1-t^{-2})(\phi_-^{(2)})\iv\phi_-^{(1)}$. {}From the conditions
on the factors $(\phi_-^{(2)})\iv$ and $\phi_-^{(1)}$ stated in 
(i) of the definition of the weak asymmetric factorization, it follows that
$\psi\in\ovl{H^1(\T)}$. Formula (\ref{f3.2}) gives
\bq
(1-t^{-2})\iv\psi(t)t^{\kappa}&=&\phi_0^{(2)}(\phi_0^{(1)})\iv,\nn
\eq
where the right hand side is an even function. Hence
\bq
(1-t^{-2})\iv\psi(t)t^{\kappa}&=&(1-t^2)\iv\psi(t\iv)t^{-\kappa}\nn
\eq
and
\bq
\psi(t)t^{2\kappa+2}&=&-\psi(t\iv).\nn
\eq
Using the fact that $\psi\in\ovl{H^1(\T)}$, we obtain $\psi=0$ if
$\kappa\le-1$ by inspecting the Fourier coefficients of $\psi$.
This is a contradiction since it would imply that $\phi_-^{(1)}=0$.
Hence $\kappa=0$. In this case it follows that $\psi(t)=\gamma(1-t^{-2})$
with $\gamma\neq0$. Hence $\phi_-^{(1)}=\gamma\phi_-^{(2)}$. 
\end{proof}

The above introduced weak asymmetric factorization is {\em not}
a kind of Wiener-Hopf factorization.
Although condition (i) on the factor $\phi_-$ implies that
$\phi_-$ can be identified with a function that is analytic and nonzero
on $\{z\in\C:|z|>1\}\cup\{\iy\}$, the factor $\phi_0$ is an even function
defined on the unit circle.

However, it is possible to relate the weak asymmetric factorization
to some kind of Wiener-Hopf factorization, which we will call weak 
antisymmetric factorization. 
In this factorization, the left and the right factor 
have a special dependence.

A function $F\in G\Li$ is said to admit a {\em weak antisymmetric factorization
in $\Lp$} if it admits a representation
\bq
F(t) &=& \phi_-(t)t^{2\kappa}\wt{\phi}_-\iv(t),\qquad t\in\T,
\eq
such that $\kappa\in\Z$ and
\begin{itemize}
\item[(i)]
$(1+t\iv)\phi_-\in\ovl{\Hp}$, $(1-t\iv)\phi_-\iv\in\ovl{\Hq}$,
\end{itemize}
where $p\iv+q\iv=1$.

Obviously, condition (i) can be rephrased by the following equivalent
condition given terms of the function $\wt{\phi}_-$:
\begin{itemize}
\item[(ii)]
$(1+t)\wt{\phi}_-\in\Hp$, $(1-t)\wt{\phi}_-\iv\in\Hq$.
\end{itemize}
Hence the function $\phi_-$ can be identified with a function that is
analytic and nonzero on $\{z\in\C:|z|>1\}\cup\{\iy\}$, while 
$\wt{\phi}_-$ can be identified with a function that is 
analytic and nonzero on $\{z\in\C:|z|<1\}$. 
In this sense, the weak antisymmetric factorization
represents a kind of Wiener-Hopf factorization.

A necessary condition for the existence of a weak antisymmetric factorization is,
of course, that $\wt{F}\iv=F$. The next result says that $\phi$ possesses
a weak asymmetric factorization if and only if the function $F=\phi\wt{\phi}\iv$
possesses a weak antisymmetric factorization.

\begin{proposition}\label{p3.2}
Let $\phi\in G\Li$ and put $F=\phi\wt{\phi}\iv$.
\begin{itemize}
\item[(a)]
If $\phi$ admits a weak asymmetric factorization,
$\phi=\phi_-t^{\kappa}\phi_0$, then the function $F$
admits a weak antisymmetric factorization with the same factor $\phi_-$ and
the same index $\kappa$.
\item[(b)]
If $F$ admits a weak antisymmetric factorization,
$F=\phi_-t^{2\kappa}\wt{\phi}_-\iv$, then $\phi$ admits a weak asymmetric
factorization with the same factor $\phi_-$, the same index $\kappa$ and
the factor $\phi_0:=t^{-\kappa}\phi_-\iv\phi$.
\end{itemize}
\end{proposition}
\begin{proof}
(a):\ Starting from the weak asymmetric factorization of $\phi$, we form
$\wt{\phi}\iv=\phi_0\iv t^{\kappa}\wt{\phi}_-\iv$ observing the fact that
$\phi_0$ is an even function, and obtain
$\phi\wt{\phi}\iv=\phi_-t^{2\kappa}\wt{\phi}_-\iv$.

(b):\ The definition of $\phi_0$ implies that $\phi=\phi_-t^{\kappa}\phi_0$.
It remains to verify that $\phi_0$ satisfies the required conditions.
First of all,
$$
\wt{\phi}_0=t^{\kappa}\wt{\phi}_-\iv\wt{\phi}=t^{\kappa}\wt{\phi}_-\iv
=t^{-\kappa}\phi_-\iv F\wt{\phi}_-\iv=t^{-\kappa}\phi_-\iv\phi=\phi_0.
$$
Hence $\phi_0$ is an even function. {}From the conditions on $\phi_-$ and
since $\phi\in G\Li$, it follows from equation
$\phi_0=t^{-\kappa}\phi_-\iv\phi$ that 
$|1-t|\phi_0\in\Lq$ and $|1+t|\phi_0\iv\in\Lp$. 
\end{proof}

{}From the last two results it follows in particular that the
index of a weak antisymmetric factorization is uniquely determined and
also the factor $\phi_-$ is uniquely determined up to a nonzero complex
constant.

%%%%%%%%%%%%%%%%%%%%%%%%%%%%%%%%%%%%%%%%%%%%%%%%%%%%%%%%%%%%%%%%%%
%%%%%%%%%%%%%%%%%%%%%%%%%%%%%%%%%%%%%%%%%%%%%%%%%%%%%%%%%%%%%%%%%%
%%%%%%%%%%%%%%%%%%%%%%%%%%%%%%%%%%%%%%%%%%%%%%%%%%%%%%%%%%%%%%%%%%

\section{The formal inverse of $\Phi(\phi)$}

We proceed with establishing some auxiliary results that we will 
need later on. Let $\cR$ stand for the linear space of all
trigonometric polynomials. Suppose that we are given a weak 
asymmetric factorization of a function $\phi\in G\Li$ with the index 
$\kappa=0$. Introduce
\bq\label{f.X1X2}
X_1=\{(1-t\iv)f(t):f\in\cR\},\quad
X_2=\{(1+t\iv)\phi_0\iv(t)f(t):f\in\cR, f(t)=f(t\iv)\}.\quad
\eq
Note that $X_1$ is dense in the Banach spaces $\Lp$, $1<p<\iy$,
and that $X_2\subseteq L^p_J(\T)$.

\begin{lemma}\label{l3.3}
Assume $\phi\in G\Li$ admits a weak asymmetric factorization in $\Lp$
with the index $\kappa=0$. Then the following assertions hold.
\begin{itemize}
\item[(a)]
The operator $B=L(\phi_0\iv)(I+J)PL(\phi_-\iv)$ is a well defined
linear operator acting from $X_1$ into $X_2$.
\item[(b)]
$\Phi(\phi)B=P|_{X_1}$.
\item[(c)]
$\ker \Phi(\phi)=\{0\}$.
\end{itemize}
\end{lemma}
\begin{proof}
(a):\
Let $f\in X_1$. We are going to compute $Bf$. First write
$f(t)=(1-t\iv)f_1(t)$ with some trigonometric polynomial $f_1$.
Hence
$$
\phi_-\iv(t)f(t)=
(1-t\iv)\phi_-\iv(t)f_1(t)
$$
i.e., a function in $\ovl{\Hq}$ times a trigonometric polynomial.
We can uniquely decompose
$$
\phi_-\iv(t)f(t)=u_1(t)+p_1(t)
$$
such that $u_1\in t\iv\ovl{\Hq}$ and $p_1$ is a polynomial.
Applying the Riesz projection to $\phi_-\iv f$ we obtain
$P(\phi_-\iv f)=p_1$. Hence
\bq\label{f3.20}
(Bf)(t) &=& \phi_0\iv(t)(p_1(t)+t\iv p_1(t\iv)).
\eq
Since $p_1(t)+t\iv p_1(t\iv)$ vanishes at $t=-1$, this expression is
$(1+t\iv)$ times a trigonometric polynomial. Now it
is easy to see that (\ref{f3.20}) is contained in $X_2$. 

(b):\ 
In order to prove (b) we need more explicit expressions for the 
terms $u_{1}$ and $p_{1}$ defined in part (a). Given
$f\in X_1$ write first
$f(t)=(1-t^{-2})f_2(t)+\alpha(1-t\iv)$ with 
some trigonometric polynomial $f_2$ and a constant $\alpha$.
Similar as before, it is possible to decompose
$$
(1-t\iv)\phi_-\iv(t)f_2(t)=u_2(t)+p_2(t)
$$
such that $u_2\in t\iv\ovl{\Hq}$ and $p_2$ is a polynomial. 
Multiplying by $(1+t^{-1})$ we obtain
$$
\phi_-\iv(t)(1-t^{-2})f_2(t)=\left((1+t\iv)u_2(t)+t\iv p_2(0)\right)+
\left((1+t\iv)p_2(t)-t\iv p_2(0)\right)
$$
Moreover, we decompose
$$
(1-t\iv)\phi_-\iv(t)=\left((1-t\iv)\phi_-\iv(t)-\phi_-\iv(\iy)\right)+
\phi_-\iv(\iy).
$$
Combining this yields
\bq
u_1(t)&=&(1+t\iv)u_2(t)+t\iv p_2(0)+
\alpha(1-t\iv)\phi_-\iv(t)-\alpha\phi_-\iv(\iy)
\nn\\
p_1(t)&=&(1+t\iv)p_2(t)-t\iv p_2(0)+\alpha\phi_-\iv(\iy).\nn
\eq
Indeed, $\phi_-\iv f=u_1+p_1$, where $u_1\in t\in\ovl{\Hq}$ and $p_1$
is a polynomial.

In order to show that $\Phi(\phi)Bf=Pf$, we have to prove
that $Pg=Pf$, where 
\bq
g(t)=\phi_-(t)(p_1(t)+t\iv p_1(t\iv)).\nn
\eq
It follows that
\bq
f(t)-g(t) &=&
\phi_-(t)u_1(t)-\phi_-(t) t\iv p_1(t\iv)
\nn\\
&=&
\phi_-(t)\left((1+t\iv)u_2(t)+t\iv p_2(0)+
\alpha(1-t\iv)\phi_-\iv(t)-\alpha\phi_-\iv(\iy)\right.\nn\\
&&\hspace{6ex}\left.
-(1+t\iv)p_2(t\iv)+p_2(0)-t\iv\alpha\phi_-\iv(\iy)
\right)
\nn\\
&=& (1+t\iv)\phi_-(t)\left(
u_2(t)+p_2(0)-p_2(t\iv)\right)-\alpha t\iv
\nn\\
&&\mbox{}+\alpha- (1+t\iv)\phi_-(t)\alpha\phi_-\iv(\iy).\nn
\eq
It can be verified easily that the expression on the
right hand side belongs to $t\iv\ovl{H^1(\T)}$.
Hence $Pf=Pg$.

(c):\
Let $f\in \ker\Phi(\phi)$. This means that $f\in L^p_J(\T)$ and
$P(\phi f)=0$. The latter can be rewritten in the form 
$f_-=\phi f$ where $f_-\in t\iv\ovl{\Hp}$. Multiplying with $\phi_-\iv$
it follows that $\phi_-\iv f_-=\phi_0f$. Moreover,
$$t(1-t\iv)\phi_-\iv(t) f_-(t)=(t-1)\phi_0(t)f(t)=:\psi(t).
$$
Since $Jf=f$ and $\phi_0$ is even, it is easily seen that $\wt{\psi}=-\psi$.
Finally, $(1-t\iv)\phi_-\iv(t)\in\ovl{\Hq}$. Hence
$\psi\in\ovl{H^1(\T)}$. It follows that $\psi=0$, thus $f_-=0$.
\end{proof}

%%%%%%%%%%%%%%%%%%%%%%%%%%%%%%%%%%%%%%%%%%%%%%%%%%%%%%%%%%%%%%%%%%
%%%%%%%%%%%%%%%%%%%%%%%%%%%%%%%%%%%%%%%%%%%%%%%%%%%%%%%%%%%%%%%%%%
%%%%%%%%%%%%%%%%%%%%%%%%%%%%%%%%%%%%%%%%%%%%%%%%%%%%%%%%%%%%%%%%%%

\section{Invertibility and asymmetric factorization}

In this section we derive necessary and sufficient conditions
for the invertibility of $M(\phi)$ in terms of an asymmetric
factorization.

\begin{lemma}\label{l4.1}
Suppose that $\Phi(\phi)$ is right invertible. Then there exist
functions $f_-\neq0$ and $f_0$ such that $f_-(t)=\phi(t)f_0(t)$, $t\in\T$, and
\be
(1+t\iv)f_-\in\ovl{\Hp},\qquad |1+t|f_0\in\Lpe.
\ee
\end{lemma}
\begin{proof}
If $\Phi(\phi)$ is right invertible, then $\im\Phi(\phi)=\Hp$. Let
$h_0\in L^p_J(\T)$ such that $\Phi(\phi)h_0=1$. Put $h_-(t)=\phi(t)h_0(t)$.
It follows that $h_-\in\ovl{\Hp}$ and $h_-\neq0$. Define
$f_-(t)=(1+t\iv)\iv h_-(t)$. Obviously, $f_-$ satisfies the required conditions
and
$$ f_-(t)=(1+t\iv)\iv \phi(t)h_0(t).$$
On defining $f_0(t)=(1+t\iv)\iv h_0(t)$, it is easily seen that
$|1+t|f_0\in\Lp$ and that $f_0$ is an even function.
\end{proof}

\begin{lemma}\label{l4.2}
Suppose that $\Psi(\psi)$ is left invertible, where $\psi(t)=\phi\iv(-t\iv)$.
Then there exist functions $g_-\neq0$ and $g_0$ such
that $g_-(t)=g_0(t)\phi\iv(t)$, $t\in\T$, and
\be
(1-t\iv)g_-\in\ovl{\Hq},\qquad |1-t|g_0\in\Lqe.
\ee
\end{lemma}
\begin{proof}
If $\Psi(\psi)$ is left invertible, then $(\Psi(\psi))^*$ is right invertible.
Identifying $(\Hp)^*$ with $\Hq$ and $(L^p_J(\T))^*$ with $L^q_J(\T)$,
it follows that
$(\Psi(\psi))^*=\Phi(\ovl{\psi})$, which is an operator acting from
$L^q_J(\T)$ into $\Hq$. It follows from the previous lemma that
there exist $f_-\neq0$ and $f_0$  such that
$$(1+t\iv)f_-\in\ovl{\Hq},\qquad |1+t|f_0\in\Lqe,$$
and 
$$ f_-(t)=\ovl{\phi\iv(-t\iv)}f_0(t).$$
Now we pass to the complex conjugate and make a substitution $t\mapsto -t\iv$.
Putting $g_-(t)=\ovl{f_-(-t\iv)}$ and $g_0(t)=\ovl{g_0(-t\iv)}$, it follows that
$$
(1-t\iv)g_-\in\ovl{\Hq},\qquad |1-t|g_0\in\Lqe,
$$ 
and $g_-(t)=\phi\iv(t)g_0(t)$, as desired.
\end{proof}

In order to present the main result, we introduce the notion of a
asymmetric factorization. A function $\phi\in G\Li$ is said to admit an
{\em asymmetric factorization in $\Lp$} if it admits a representation
\bq\label{f3.asf}
\phi(t)&=& \phi_-(t)t^{\kappa}\phi_0(t),\qquad t\in\T,
\eq
such that $\kappa\in\Z$ and
\begin{itemize}
\item[(i)]
$(1+t\iv)\phi_-\in\ovl{\Hp}$, $(1-t\iv)\phi_-\iv\in\ovl{\Hq}$,
\item[(ii)]
$|1-t|\phi_0\in\Lqe$, $|1+t|\phi_0\iv\in\Lpe$,
\item[(iii)]
the linear operator $B=L(\phi_0\iv)(I+J)PL(\phi\iv_-)$ acting
from $X_1$ into $X_2$ extends to a linear bounded operator $\wt{B}$
acting from $\Lp$ into $L^p_J(\T)$.
\end{itemize}
Here $p\iv+q\iv=1$, and $X_1$ and $X_2$ are the spaces (\ref{f.X1X2}).

An equivalent formulation of condition (iii) is the following:
\begin{itemize}
\item[(iii*)]
there exists a constant $M$ such that
$\|Bf\|_{L^p_J(\T)}\le M\|f\|_{\Lp}$ for all $f\in X_1$.
\end{itemize}
Notice  here that $X_1$ is dense in $\Lp$ and $X_2\subseteq L^p_J(\T)$.

\begin{theorem}\label{t4.3}
Let $\phi\in G\Li$. The operator $M(\phi)$ is invertible on $\Hp$ if and
only if $\phi$ admits an asymmetric factorization in $\Lp$ with the index 
$\kappa=0$.
\end{theorem}
\begin{proof}
$\Rightarrow$:\
If $M(\phi)$ is invertible, then by Proposition \ref{p2.4} the operators
$\Phi(\phi)$ and $\Psi(\psi)$ are also invertible. Applying Lemma \ref{l4.1} and
Lemma \ref{l4.2} it follows that $f_-=\phi f_0$ and $g_-=g_0\phi\iv$
with the properties stated there. Combining this yields
$g_-f_-=g_0f_0$. Now we are in a similar situation as in the proof of
Proposition \ref{p3.1} (see also (\ref{f3.16})). In the very same way one can prove
that $g_-f_-=g_0f_0=:\gamma$ is a constant. It must be nonzero since otherwise
$g_-=0$ or $f_-=0$ (cf.~the F.~and M.~Riesz Theorem).
Now we put $\phi_-=f_-=\gamma g_-\iv$ and $\phi_0=f_0\iv=g_0\gamma\iv$.
{}From the conditions on the functions in Lemma \ref{l4.1} and
Lemma \ref{l4.2}, the conditions (i) and (ii) in the above definition
of the asymmetric factorization follow. Hence we have shown that
$\phi$ admits a weak asymmetric factorization with $\kappa=0$.

{}From Lemma \ref{l3.3}(a) it follows that the operator $B$ is well defined.
Assertion (b) of the same lemma implies that (since $\Phi(\phi)$ is invertible)
$$B\;=\;\Phi(\phi)\iv P|_{X_1}.$$
Both the left and right hand side in this equation represent operators defined
on $X_1$, which is dense in $\Lp$. Obviously, the right hand side can
be extended by continuity to a linear bounded operator acting
from $\Lp$ into $\Hp$ (since it is the restriction of such an operator).
Hence so can be the operator $B$.

$\Leftarrow$:\
Assume that conditions (i)-(iii) are satisfied. Denote the continuous
extension of the operator $B$ by $\wt{B}$. We can apply Lemma \ref{l3.3}(b).
Since $X_1$ is dense in $\Lp$ it follows that
$$\Phi(\phi) \wt{B} = P,$$
where the operators are defined on $\Lp$. Hence $\Phi(\phi)$ is right invertible,
the right inverse being $\wt{B}|_{\Hp}$ (the restriction of $\wt{B}$ onto $\Hp$).
Now we multiply the above equation from the right with $\Phi(\phi)$ and obtain
$\Phi(\phi)\wt{B}|_{\Hp}\Phi(\phi)=\Phi(\phi)$, i.e.,
$$\Phi(\phi)\left(\wt{B}|_{\Hp}\Phi(\phi)-I\right)=0.$$
{}From the triviality of the kernel of $\Phi(\phi)$ (see Lemma \ref{l3.3}(c)),
we obtain $\wt{B}|_{\Hp}\Phi(\phi)=I$. Hence $\Phi(\phi)$ is invertible and its
inverse is just $\wt{B}|_{\Hp}$.
\end{proof}

%%%%%%%%%%%%%%%%%%%%%%%%%%%%%%%%%%%%%%%%%%%%%%%%%%%%%%%%%%%%%%%%%%
%%%%%%%%%%%%%%%%%%%%%%%%%%%%%%%%%%%%%%%%%%%%%%%%%%%%%%%%%%%%%%%%%%
%%%%%%%%%%%%%%%%%%%%%%%%%%%%%%%%%%%%%%%%%%%%%%%%%%%%%%%%%%%%%%%%%%

\section{Fredholm theory and asymmetric factorization}

In order to establish the Fredholm theory in terms of the above introduced
asymmetric factorization, it is necessary to prove an assertion about
the kernel and the cokernel of $M(\phi)$.
It is the analogue of Coburn's result for Toeplitz operators.

We begin with the following observation.
Given $\phi\in L^\iy(\T)$, define
\bq
K&=&\left\{\,t\in\T\;\left|\;\phi(t)=\phi(t\iv)=0\,\right.\right\}.
\eq
This definition of course depends on the choice of representatives 
for $\phi,$ however the following remarks are independent of that 
choice.
Because the characteristic function $\chi_K$ is an even function,
it follows from (\ref{f.M}) that $M(\chi_K)$ is a projection.
Moreover, from the same relation 
we obtain that $M(\phi)M(\chi_K)=M(\phi\chi_K)=0$ and thus
\bq
\im M(\chi_K)&\subseteq&\ker M(\phi).
\eq
If the set $K$ has Lebesgue measure zero, then obviously 
$M(\chi_K)=0$,
whereas if $K$ has a positive Lebesgue measure, then the image of 
$M(\chi_K)$ is
infinite dimensional. The latter fact can be seen by decomposing $K$ 
into
pairwise disjoint and even sets $K_1,\dots,K_n$ with positive 
Lebesgue measure.
Then $M(\chi_K)=M(\chi_{K_1})+\dots+M(\chi_{K_n})$, where
$M(\chi_{K_1}),\dots,M(\chi_{K_n})$ are mutually orthogonal projections.
By Proposition \ref{p2.1} all these projections are nonzero.
Hence $\dim\im M(\chi_K)\ge n$ for all $n$, i.e., $\dim\im M(\chi_K)=\iy$.

\begin{proposition}\label{p2.5}
Let $\phi\in L^\iy(\T)$ and let $K$ be as above. Then
$\ker M(\phi)=\im M(\chi_K)$ or $\ker M^*(\phi)=\{0\}$.
\end{proposition}
\begin{proof}
Obviously, $M(\phi)$ and its adjoint can be written in the form:
$$ M(\phi)\;=\; PL(\phi)(I+J)P,\qquad\qquad
M^*(\phi)\;=\; P(I+J)L(\overline{\phi})P.$$
Suppose that we are given functions $f_+\in H^p(\T)$ and $g_+\in H^q(\T)$
such that  $M(\phi)f_+= 0$
and $M^*(\phi)g_+= 0$ with $g_+\neq0$. We have to show that 
$f_+\in\im M(\chi_K)$.
Introducing the functions
$$\ba{rclcrcl}
f(t) &=& f_+(t)+t\iv f_+(t\iv), &\qquad& f_-(t) &=& 
\phi(t)f(t),\\[1ex]
g(t) &=& \overline{\phi}(t)g_+(t), && g_-(t) &=& g(t)+t\iv g(t\iv),
\ea$$
it follows that $f_-\in t\iv\ovl{\Hp}$ and $g_-\in t\iv\ovl{\Hq}$.
{}From the definition of $g_-$ we obtain that
$t\iv g_-(t\iv) = g_-(t)$. Hence $g_-=0$ by checking the Fourier 
coefficients.
It follows that $t\iv g(t\iv )=-g(t)$. On the other hand, the 
definition of $f$
says that $t\iv f(t\iv )=f(t)$. This implies that
$(f\overline{g})(t\iv)=-(f\overline{g})(t)$.
Also from the above relations we conclude that
$f\overline{g}=f\phi\overline{g}_+= f_-\overline{g}_+$.
Because $f_-\in t\iv\ovl{\Hp}$ and $\overline{g}_+\in 
\overline{H^q(\T)}$, we
have $f_-\overline{g}_+\in t\iv\ovl{H^1(\T)}$.  Considering again the Fourier
coefficients, this shows $f\overline{g}=f_-\overline{g}_+=0$. 
Since $g_+\in H^q(\T)$ and $g_+\neq0$ the F. and M. Riesz Theorem
says that $g_+(t)\neq 0$ almost everywhere on $\T$, and thus we 
obtain $f_-=0$.
Hence $\phi f=0$. Because $t\iv f(t\iv )=f(t)$, we have also
$\wt{\phi}f=0$. The definition of the set $K$ now implies that
$(1-\chi_K)f=0$. Noting that $M(1-\chi_K)f_+=0$, we finally arrive at
$f_+\in\im M(\chi_K)$.
\end{proof}

Combining the previous proposition with Proposition \ref{p2.4}
we obtain the following result.
\begin{corollary}\label{c5.2}
Let $\phi\in L^\iy(\T)$. If $M(\phi)$ is Fredholm on $\Hp$, then $M(\phi)$ has 
a trivial kernel or a trivial cokernel.
\end{corollary}

{}From this corollary we can obtain another interesting result, pertaining
to the relation between invertibility and Fredholmness of $M(\phi)$.

\begin{corollary}\label{c5.3}
Let $\phi\in L^\iy(\T)$. Then $M(\phi)$ is invertible on $\Hp$ if and only if
$M(\phi)$ is Fredholm on $\Hp$ and has index zero.
\end{corollary}

The desired Fredholm criteria in terms of the asymmetric factorization
is established next.

\begin{theorem}\label{t5.4}
Let $\phi\in G\Li$. The operator $M(\phi)$ is a Fredholm operator on $\Hp$
if and only if the function $\phi$ admits an asymmetric factorization 
{\rm (\ref{f3.asf})} in $\Lp$. In this case the defect numbers are given by
\be
\dim\ker M(\phi)=\max\{0,-\kappa\},\qquad
\dim\ker (M(\phi))^*=\max\{0,\kappa\}.
\ee
\end{theorem}
\begin{proof}
$\Rightarrow$:\ Assume that $M(\phi)$ is Fredholm with index $-\kappa$.
We consider the function $\psi(t)=t^{-\kappa}\phi(t)$. Using the fact that a
Hankel operator with a continuous generating function is compact and the
identity (\ref{f.Mx}), it follows that
$$M(\psi)=M(t^{-\kappa})M(\phi)+\mbox{compact}.$$
Taking into account that $M(t^{-\kappa})=T(t^{-\kappa})+\mbox{compact}$
is Fredholm with index $\kappa$, it follows that $M(\psi)$
is a Fredholm operator with index zero. Now we are in a position to apply
Corollary \ref{c5.3} in order to see that $M(\psi)$ is invertible.
Theorem \ref{t4.3} implies that $\psi$ possesses an asymmetric factorization
with index zero. Hence $\phi$ possesses an asymmetric factorization with
index $\kappa$.

$\Leftarrow$:\ Assume that $\phi$ possesses an asymmetric factorization with
index $\kappa$. Then we consider again the function $\psi(t)=t^{-\kappa}\phi(t)$,
which thus possesses an asymmetric factorization with index zero. Theorem 
\ref{t4.3} implies that $M(\psi)$ is invertible. Very similar as above it follows
that $M(\phi)$ is Fredholm and has Fredholm index $-\kappa$.

As to the defect numbers, we known that $M(\phi)$ has the Fredholm index 
$-\kappa$. On the other hand (by Corollary \ref{c5.2}), the kernel or the
cokernel of $M(\phi)$ are trivial. This implies the formulas for the
defect numbers.
\end{proof}

%%%%%%%%%%%%%%%%%%%%%%%%%%%%%%%%%%%%%%%%%%%%%%%%%%%%%%%%%%%%%%%%%
%%%%%%%%%%%%%%%%%%%%%%%%%%%%%%%%%%%%%%%%%%%%%%%%%%%%%%%%%%%%%%%%%
%%%%%%%%%%%%%%%%%%%%%%%%%%%%%%%%%%%%%%%%%%%%%%%%%%%%%%%%%%%%%%%%%

\section{Fredholm and invertibility theory 
for piecewise continuous functions}

The goal of this section is to find practical criteria that determine 
whether or not $M(\phi)$ is an invertible operator on $\Hp$ when $\phi$ 
is a piecewise continuous function. 

The first part of this section is devoted to the Fredholmness of
$M(\phi)$ for piecewise continuous functions $\phi$.
The results are immediate consequences from \cite{RS}.
Let us start with the following definitions.
The Mellin transform $\cM:L^p(\Rp)\to L^p(\R)$ is defined by
\bq
(\cM f)(x) &=& \int_0^\iy \xi^{-1+1/p-ix}f(\xi)\,d\xi,\quad x\in\R.
\eq
For $\phi\in L^\iy(\R)$ the Mellin convolution operator
$M^0(\phi)\in\cL(L^p(\Rp))$ is given by
\bq
M^0(\phi)f &=& \cM\iv(\phi(\cM f)).
\eq
Let $S$ be the singular integral operator acting on the space $L^p(\Rp)$,
\bq
(Sf)(x) &=& \frac{1}{\pi i}\int_0^\iy \frac{f(y)}{y-x}\,dy,\quad x\in\Rp,
\eq
where the singular integral has to be understood as the Cauchy principal value,
and let $N$ stand for the integral operator acting on $L^p(\Rp)$ by
\bq
(Nf)(x) &=& \frac{1}{\pi i}\int_0^\iy \frac{f(y)}{y+x}\,dy,\quad x\in\Rp.
\eq
It is well known that $S$ and $N$ can be expressed as Mellin convolution
operators
\be
S=M^0(s), \qquad N=M^0(n)
\ee
with the generating functions
\be
s(z)=\coth((z+i/p)\pi),\quad
n(z)=\sinh\iv((z+i/p)\pi),\quad z\in\R
\ee
(see, e.g., Section 2 of \cite{RS}). Notice that $s$ and $n$ are continuous
on $\R$ and possess the limits $s(\pm\iy)=\pm1$ and $n(\pm\iy)=0$ for
$z\to\pm\iy$. Moreover, $s^2-n^2=1$.

Let $S_\T$ stand for the singular integral operator acting on $L^p(\T)$,
\bq
(S_\T f)(t) &=& \frac{1}{\pi i} \int_\T \frac{f(s)}{s-t}ds, \quad t\in\T,
\eq
and let $W$ stand for the ``flip'' operator 
\bq
(W f)(t) &=& f(t\iv)
\eq
acting on $L^p(\T)$. Finally, for a piecewise continuous function $\phi\in PC$
defined on the unit circle denote by
\bq
\phi_\pm(t) &=& \lim\limits_{\varepsilon\to\pm0}\phi(te^{i\varepsilon})
\eq
the one-sided limits at a point $t\in \T$. 
As before let $L(\phi)$ stand for the multiplication operator
on $L^p(\T)$ with the generating function $\phi$.

The smallest closed subalgebra of $\cL(L^p(\T))$, which contains the operators
$S$, $W$ and $L(\phi)$ for $\phi\in PC$ will be denoted by
$\cS^p(PC)$.
Specializing Theorem 9.1 of \cite{RS} to the situation where we are 
interested in, we obtain the following result. Therein we
put $\T_+=\{\tau\in\T:\Im \tau>0\}$.

\begin{theorem}\label{t7.1}\
\begin{itemize}
\item[(a)]
For $\tau\in\{-1,1\}$, there exists a homomorphism 
$H_\tau:\cS^p(PC)\to \cL((L^p(\T))^2)$, which acts on the generating 
elements as follows:
\bq
H_\tau (S_\T) &=& 
\left(\ba{cc} S&-N\\ N&-S\ea\right),\qquad
H_\tau (W) \;\;=\;\;
\left(\ba{cc} 0&I\\ I&0\ea\right),
\nn\\[2ex]
H_\tau (L(\phi))  &=& 
\diag(\phi_+(\tau)I,\phi_-(\tau)I).\nn
\eq
\item[(b)]
For $\tau\in\T_+$, there exists a homomorphism 
$H_\tau:\cS^p(PC)\to \cL((L^p(\T))^4)$, which acts on the generating 
elements as follows:
\bq
H_\tau (S_\T) &=& 
\left(\ba{cccc} S&-N&0&0\\ N&-S&0&0\\ 0&0&S&-N\\0&0&N&-S\ea\right),\qquad
H_\tau (W) \;\;=\;\;
\left(\ba{cccc} 0&0&0&I\\ 0&0&I&0\\ 0&I&0&0\\ I&0&0&0\ea\right),
\nn\\[2ex]
H_\tau (L(\phi))  &=& 
\diag(\phi_+(\tau)I,\phi_-(\tau)I,\phi_+(\bar{\tau})I,\phi_-(\bar{\tau})I).\nn
\eq
\item[(c)]
An operator $A\in \cS^p(PC)$ is Fredholm if and only if the operators
$H_\tau (A)$ are invertible for all $\tau\in\{-1,1\}\cup\T_+$.
\end{itemize}
\end{theorem}

Now we apply this theorem in order to study the Fredholmness of
the operator $M(\phi)$.

\begin{theorem}\label{t7.2new}
Let $\phi\in PC$. Then $M(\phi)$ is Fredholm on $\Hp$ if and only if
$\phi_\pm(\tau)\neq0$ for all $\tau\in\T$ and
and if the following conditions are satisfied:
\bq
\frac{1}{2\pi}\arg\left(\frac{\phi_-(1)}{\phi_+(1)}\right) &\notin& 
\frac{1}{2p}+\Z,\label{f.38}\\
\frac{1}{2\pi}\arg\left(\frac{\phi_-(-1)}{\phi_+(-1)}\right) &\notin& 
\frac{1}{2}+\frac{1}{2p}+\Z,\label{f.39}\\
\frac{1}{2\pi}\arg\left(
\frac{\phi_-(\tau)\phi_-(\bar{\tau})}{\phi_+(\tau)\phi_+(\bar{\tau})}\right) 
&\notin& \frac{1}{p}+\Z\qquad
\mbox{ for each }\tau\in\T_+.\label{f.40}
\eq
\end{theorem}
\begin{proof}
Observe that $M(\phi)=PL(\phi)(I+J)P$ is defined on $\Hp$.
Moreover, as has been shown in the proof of Proposition \ref{p2.4},
the operators $(I+J)P:\Hp\to L^p_{J}(\T)$ and 
$\frac{1}{2}P(I+J):L^p_{J}(\T)\to\Hp$
are inverse to each other. Hence $M(\phi)$ is Fredholm on $\Hp$ if and only if
the operator $A=(I+J)PM(\phi)P(I+J)$ acting on $L^p_{J}(\T)$ 
is Fredholm. Taking into account that $PJ=J(I-P)$ and $J^2=I$, 
a simple computation gives
\bq
A &=& (I+J)PL(\phi)(I+J).\nn
\eq
Notice that $P=(I+S_{\T})/2$ and $J=L(t\iv)W$. The space $L^p_{J}(\T)$
is by definition the image of the operator $(I+J)$. The operator $A$ belongs
to $\cS^p(PC)$ and thus we can apply Theorem \ref{t7.1} in order to examine
the Fredholmness of the operator $A$. The results is that
$A$ is Fredholm on $L^p_{J}(\T)$ if and only if for each $\tau\in\T_+\cup
\{-1,1\}$ the operators $H_\tau(A)$ are invertible on image of the operators
$H_\tau(I+J)$.

For $\tau\in\{-1,1\}$ we obtain that
\bq
H_\tau(I+J) &=&
\left(\ba{cc} I&\tau I\\\tau I&I\ea\right)
\;\;=\;\;
\left(\ba{c}I\\\tau I\ea\right) \left(\ba{cc} I&\tau I\ea\right)\nn
\eq
and 
\bq
H_\tau(A)  &=&\frac{1}{2}\left(\ba{cc} I&\tau I\\\tau I&I\ea\right)
\left(\ba{cc} (I+S)\phi_+(\tau) & -N\phi_-(\tau)\\ 
N\phi_+(\tau)&(I-S)\phi_-(\tau)\ea\right)
\left(\ba{cc} I&\tau I\\\tau I&I\ea\right).\nn
\eq
Thus $H_\tau(A)$ is invertible on the image of $H_\tau(I+J)$ 
if and only if the operator
\bq
B_\tau &=& \left(\ba{cc} I&\tau I\ea\right)
\left(\ba{cc} (I+S)\phi_+(\tau) & -N\phi_-(\tau)\\ 
N\phi_+(\tau)&(I-S)\phi_-(\tau)\ea\right)
\left(\ba{c}I\\\tau I\ea\right)\nn
\eq
is invertible on $L^p(\T)$. Obviously,
\bq
B_\tau &=& (I+S)\phi_+(\tau)+(I-S)\phi_-(\tau)
+\tau N(\phi_+(\tau)-\phi_-(\tau)),\nn
\eq
which is a Mellin convolution operator $B_\tau =M^0(b_\tau)$ with the generating
function
\bq
b_\tau(z) &=& (1+s(z)+\tau n(z))\phi_+(\tau)+(1-s(z)-\tau n(z))\phi_-(\tau),
\qquad z\in\R.\nn
\eq
Thus $B_\tau$ is invertible if and only if $b_\tau(\pm\iy)\neq0$ and
$b_\tau(z)\neq0$ for all $z\in\R$. Obviously, $b_\tau(\pm\iy)=2\phi_\pm(\tau)$.
The second condition is equivalent to
\bq
\frac{\phi_-(\tau)}{\phi_+(\tau)} &\neq& 
\frac{s(z)+\tau n(z)+1}{s(z)+\tau n(z)-1}
\;\;=\;\;\frac{e^{(z+i/p)\pi}+\tau}{e^{-(z+i/p)\pi}+\tau}
\;\;=\;\; \tau e^{(z+i/p)\pi}.\nn
\eq
This implies condition (\ref{f.38}) and (\ref{f.39}).

For $\tau\in\T_+$ we have
\bq
H_\tau(I+J) &=&
\underbrace{
  \left(\ba{cccc} I&0&0&\bar{\tau}I\\0&I&\bar{\tau}I&0\\0&\tau I&I&0\\
  \tau I&0&0&I\ea\right) }_{\textstyle \mathbf{X}}
\;\;=\;\;
\underbrace{
  \left(\ba{cc} I&0\\0&I\\0&\tau I\\\tau I&0\ea
  \right)}_{\textstyle \mathbf{Y}}
\underbrace{
  \left(\ba{cccc}I&0&0&\bar{\tau}I\\0&I&\bar{\tau}I&0\ea
  \right)}_{\textstyle \mathbf{Z}}
\nn 
\eq
and 
\bq
&& 
H_\tau(A) \;\;=\;\;
\frac{1}{2} \mathbf{X}
\left(\ba{cccc} (I+S)\phi_+(\tau) & -N\phi_-(\tau) &0&\\ 
N\phi_+(\tau)&(I-S)\phi_-(\tau)&0&\\ 
0&0&(I+S)\phi_+(\bar{\tau}) & -N\phi_-(\bar{\tau})\\ 
0&0&N\phi_+(\bar{\tau})&(I-S)\phi_-(\bar{\tau})\ea\right)
\mathbf{X}.\nn
\eq
Thus $H_\tau(A)$ is invertible on the image of $H_\tau(I+J)$ 
if and only if the operator $B_\tau$ given by
$$
\mathbf{Z}
\left(\ba{cccc} (I+S)\phi_+(\tau) & -N\phi_-(\tau) &0&\\ 
N\phi_+(\tau)&(I-S)\phi_-(\tau)&0&\\ 
0&0&(I+S)\phi_+(\bar{\tau}) & -N\phi_-(\bar{\tau})\\ 
0&0&N\phi_+(\bar{\tau})&(I-S)\phi_-(\bar{\tau})\ea\right)
\mathbf{Y}
$$
is invertible. It follows that
\bq
B_\tau &=& 
\left(\ba{cc}(I+S)\phi_+(\tau) & -N\phi_-(\tau)\\ 
N\phi_+(\tau)&(I-S)\phi_-(\tau)\ea\right)+
\left(\ba{cc}(I-S)\phi_-(\bar{\tau})&N\phi_+(\bar{\tau})\\
-N\phi_-(\bar{\tau})&(I+S)\phi_+(\bar{\tau})\ea\right).\nn
\eq
This operator is a Mellin convolution operator with the generating function
\bq
b_\tau(z) &=& 
\left(\ba{cc}(1+s(z))\phi_+(\tau)+(1-s(z))\phi_-(\bar{\tau})& 
n(z)(\phi_+(\bar{\tau})-\phi_-(\tau))\\ 
n(z)(\phi_+(\tau)-\phi_-(\bar{\tau}))&
(1-s(z))\phi_-(\tau)+(1+s(z))\phi_+(\bar{\tau})\ea\right).\nn
\eq
Note that
$$
b_\tau(+\iy) = 2\left(\ba{cc}\phi_+(\tau)&0\\0&\phi_+(\bar{\tau})\ea\right),
\qquad
b_\tau(-\iy) = 2\left(\ba{cc}\phi_-(\bar{\tau})&0\\0&\phi_-(\tau)\ea\right).
$$
Finally,
\bq
\det b_\tau &=&
(1+s)^2\phi_+(\tau)\phi_+(\bar{\tau})+(1-s)^2\phi_-(\tau)\phi_-(\bar{\tau})
+(1-s^2)(\phi_+(\tau)\phi_-(\tau)+\phi_+(\bar{\tau})\phi_-(\bar{\tau}))
\nn\\
&&
-n^2(\phi_+(\tau)-\phi_-(\bar{\tau}))(\phi_+(\bar{\tau})-\phi_-(\tau))
\nn\\
&=& (1+s)^2\phi_+(\tau)\phi_+(\bar{\tau})+(1-s)^2\phi_-(\tau)\phi_-(\bar{\tau})
+(1-s^2)(\phi_+(\tau)\phi_+(\bar{\tau})+\phi_-(\tau)\phi_-(\bar{\tau}))
\nn\\
&=& 2(1+s)\phi_+(\tau)\phi_+(\bar{\tau})+2(1-s)\phi_-(\tau)\phi_-(\bar{\tau}).\nn
\eq
This expression is nonzero if and only if
\bq
\frac{\phi_-(\tau)\phi_-(\bar{\tau})}{\phi_+(\tau)\phi_+(\bar{\tau})}
&\neq&
\frac{s+1}{s-1}\;\;=\;\;e^{2(z+i/p)\pi}.\nn
\eq
{}From this condition (\ref{f.40}) follows.
\end{proof}

In the rest of this section we are going to establish the criteria for
an operator $M(\phi)$ to be invertible on $\Hp$ in the case where
$\phi$ is a piecewise continuous function with a finite number jumps.
The proof is based partly on the Fredholm criteria that we have just established
and partly on the notion of a weak asymmetric factorization.

It is well known that any piecewise continuous and
nonvanishing function with a finite number of discontinuities at the
points $t_{1}=e^{i\theta_{1}},\dots,t _{R}=e^{\theta_{R}}$ 
can be written as a product
\be
\phi(e^{i\theta}) = b(e^{i\theta})\prod_{r=1}^{R}t_{\beta_{r}}(e^{i(\theta 
-\theta_{r}})) 
\ee
where $b$ is a nonvanishing continuous function and 
\be
t_{\beta}(e^{i\theta}) = \exp(i\beta(\theta-\pi)), 
\qquad 0 < \theta < 2\pi.
\ee
The parameters in
the formula are useful to decide invertibility.  For example, it is 
well known \cite{BS2}, that the Toeplitz operator $T(\phi)$ is invertible 
on $\Hp$ if and only if it
can be represented in the above form with parameters satisfying
\be
-1/q<\Re \beta _{r}<1/p
\ee
for all $1\le r\le R$ and a continuous nonvanishing function $b$
with winding number zero.  For the corresponding result for
the operators $M(\phi)$ we prepare with the following proposition.

\begin{proposition}\label{p7.3}
Let $\psi$ be a function of the form
\be\label{f.psi}
\psi(e^{i\theta}) = 
t_{\beta^{+}}(e^{i\theta})
t_{\beta^{-}}(e^{i(\theta-\pi)})\prod_{r=1}^{R}
t_{\beta_{r}^{+}}(e^{i(\theta -\theta_{r})}) 
t_{\beta_{r}^{-}}(e^{i(\theta +\theta_{r})}) 
\ee
where $\theta_1,\dots,\theta_R\in(0,\pi)$ are distinct and
\begin{enumerate}
\item[(i)]
$-1/q < \Re(\beta_{r}^{+}+\beta_{r}^{-}) < 1/p $ 
\item[(ii)] 
$ -1/2 -1/2q < \Re\beta^{+} < 1/2p $ and 
$-1/2q < \Re\beta^{-} < 1/2 +1/2p .$
\end{enumerate}
Then $\psi$ admits a weak asymmetric factorization in $L^p(\T)$ with the
index $\kappa=0$.
\end{proposition}
\begin{proof}
First observe that
\bq
t_{\beta,\theta_{r}}(t) &=&
\left(1-tt_{r}\iv\right)^\beta
\left(1-t_{r}t\iv\right)^{-\beta}\nn\\
&=&
\left[\left(1-t_{r}t\iv\right)^{-\beta}
\left(1-t\iv t_{r}\iv\right)^{-\beta}\right]
\left[\left(1-tt_{r}\iv\right)^\beta
\left(1-t\iv t_{r}\iv\right)^\beta\right].\nn
\eq
Hence we can factor $\psi=\psi_-\psi_0$ with
\bq
\psi_-&=&
\left(1-t\iv\right)^{-2\beta^+}
\left(1+t\iv\right)^{-2\beta^-}\nn\\
&&\prod_{r=1}^R
\left(1-t_{r}t\iv\right)^{-\beta_{r}^+-\beta_{r}^-}
\left(1-t\iv t_{r}\iv\right)^{-\beta_{r}^+-\beta_{r}^-},\nn
\\
\psi_0&=&
|1-t|^{2\beta^+}|1+t|^{2\beta^-}\nn\\
&&\prod_{r=1}^R
\left(1-tt_{r}\iv\right)^{\beta_{r}^+}
\left(1-t\iv t_{r}\iv\right)^{\beta_{r}^+}
\left(1-tt_{r}\right)^{\beta_{r}^-}
\left(1-t_{r}t\iv\right)^{\beta_{r}^-}.\nn
\eq
Since $-2\Re\beta^+>-1/p$, $1-2\Re\beta^->-1/p$ and 
$-\Re\beta_{r}^+-\Re\beta_{r}^->-1/p$ we have
$$
(1+t\iv)\psi_-\in\ovl{\Hp},
$$
and since $1+2\Re\beta^+>-1/q$, $2\Re\beta^->-1/q$ and 
$\Re\beta_{r}^++\Re\beta_{r}^->-1/q$ we have
$$
(1-t\iv)\psi_-\iv\in\ovl{\Hq}.
$$
Obviously, $\psi_0$ is even and thus it easily follows that
$|1-t|\phi_0\in\Lqe$ and $|1+t|\phi_0\iv\in\Lpe$.
\end{proof}

Now we are prepared to established the invertibility
criteria for the operators $M(\phi)$.

\begin{theorem} \label{t7.2}
Suppose that $\phi$ has finitely many jump discontinuities. Then
$M(\phi)$ is invertible on $\Hp$ if and only if $\phi$ can be written in the form
\be\label{f.phi}
\phi(e^{i\theta}) = 
b(e^{i\theta})t_{\beta^{+}}(e^{i\theta})
t_{\beta^{-}}(e^{i(\theta-\pi)})\prod_{r=1}^{R}
t_{\beta_{r}^{+}}(e^{i(\theta -\theta_{r})}) 
t_{\beta_{r}^{-}}(e^{i(\theta +\theta_{r})}) 
\ee
where $b$ is a continuous nonvanishing function with winding number
zero, $\theta_{1}, \dots, \theta_{R} \in (0, \pi)$ are distinct, and
\begin{enumerate}
\item[(i)]
$-1/q < \Re(\beta_{r}^{+}+\beta_{r}^{-}) < 1/p $ 
\item[(ii)] 
$ -1/2 -1/2q < \Re\beta^{+} < 1/2p $ and 
$-1/2q < \Re\beta^{-} < 1/2 +1/2p .$
\end{enumerate}
\end{theorem}
\begin{proof}
We begin with the ``if'' part of the proof.  Assume that the parameters 
satisfy the above inequalities, let $b$ be a nonvanishing continuous function
with winding number zero, and let $\psi$ be the function
(\ref{f.psi}). Observe that $\phi=b\psi$.

It follows from Theorem \ref{t7.2new} that $M(\psi)$ is Fredholm.
Hence, by Theorem \ref{t5.4}, the function $\psi$ admits an asymmetric
factorization in $L^p(\T)$. Moreover, we know from 
Proposition \ref{p7.3} that $\psi$ admits a weak asymmetric
factorization in $L^p(\T)$ with $\kappa=0$. Due to the uniqueness
of (weak) asymmetric factorizations in the sense of Proposition \ref{p3.1},
it follows that $\psi$ admits an asymmetric 
factorization in $L^p(\T)$ with $\kappa=0$. Thus, again by Theorem \ref{t5.4},
$M(\psi)$ is Fredholm with index zero. 

It is well know that under the above assumptions on the function $b$
the Toeplitz operator $T(b)$ is Fredholm with index zero and 
the Hankel operator $H(b)$ is compact. Hence $M(b)=T(b)+H(b)$ is also
a Fredholm operator with index zero.
Taking (\ref{f.Mx}) into account, we have
\bq\label{f.50}
M(\phi) &=& M(b)M(\psi)+H(b)M(\wt{\psi}-\psi),
\eq
where the last term is compact. Thus $M(\phi)$ is Fredholm and has index
zero. Corollary \ref{c5.3} now implies that $M(\phi)$ is invertible.

In order to justify the other direction we first note that if $M(\phi)$ is 
invertible then $\phi$ satisfies the conditions 
(\ref{f.38})--(\ref{f.40}) of Theorem \ref{t7.2new}.
Hence it is possible to choose parameters $\beta^{+}, \beta^{-},
\beta_{r}^{+}, \beta_{r}^{-}$ in such a way that the fulfill the conditions
(i) and (ii) and such that
$$
\frac{\phi_-(1)}{\phi_+(1)} = \exp(2\pi i\beta^{+}),\qquad
\frac{\phi_-(-1)}{\phi_+(-1)} = \exp(2\pi i\beta^{-}),
$$
$$
\frac{\phi_-(t_r)}{\phi_+(t_r)} = \exp(2\pi i\beta_r^{+}),\qquad
\frac{\phi_-(t_r\iv)}{\phi_+(t_r\iv)} = \exp(2\pi i\beta_r^{-}).
$$
With those parameters we can represent $\phi$ in the form (\ref{f.phi}),
where it follows that $b$ is a continuous nonvanishing function.
It remains to show that the winding number of $b$ is zero.

By what we have proved in the ``if'' part of the theorem, we know that
for the function $\psi$ defined by (\ref{f.psi}) the operator
$M(\psi)$ is invertible. 
We rely on (\ref{f.50}) and can conclude that the Fredholm index of
$M(b)$ is zero. Thus $b$ has winding number zero.
\end{proof}

The authors have studied other equivalent conditions for invertibility of 
$M(\phi)$ which generalize the idea of $A_{p}$-conditions for 
weighted spaces. These results are forth-coming and will be contained 
in the sequel to this paper.

%%%%%%%%%%%%%%%%%%%%%%%%%%%%%%%%%%%%%%%%%%%%%%%%%%%%%%%%%%%%%%%%%%
%%%%%%%%%%%%%%%%%%%%%%%%%%%%%%%%%%%%%%%%%%%%%%%%%%%%%%%%%%%%%%%%%%
%%%%%%%%%%%%%%%%%%%%%%%%%%%%%%%%%%%%%%%%%%%%%%%%%%%%%%%%%%%%%%%%%%

\end{document}